\documentclass[twoside,a4paper,12pt,centertags,reqno]{amsart} 
\usepackage{amsmath,amssymb,verbatim,vmargin}
\usepackage{color}
\usepackage{tikz}
\usepackage{tikz-cd}
\usetikzlibrary{matrix,calc}
\usepackage{hyperref}
\hypersetup{colorlinks,bookmarks=true,linktocpage=true,citecolor=blue,linkcolor=magenta}

\usepackage{eulervm}   

\allowdisplaybreaks 
\usepackage{color}
\usepackage[all]{xy} 

\usepackage{mathtools}
\usepackage{MnSymbol} 

\theoremstyle{plain}
\newtheorem{thm}{Theorem}[section]
\newtheorem{lem}[thm]{Lemma}

\theoremstyle{definition}

\newtheorem{rem}[thm]{Remark}

\usepackage{enumerate}
\usepackage{enumitem}

\newcommand{\bRn}{\mathbb{R}^n}

\newcommand{\pd}{\partial}

\newcommand{\bC}{{\mathbb C}}

\newcommand{\bR}{{\mathbb R}}

\newcommand{\bL}{{\mathbb L}}

\def\barint_#1{\mathchoice
            {\mathop{\vrule width 6pt
height 3 pt depth -2.5pt
                    \kern -9.5pt
\intop \kern -4pt}\nolimits_{#1}}%
            {\mathop{\vrule width 5pt height
3 pt depth -2.6pt
                    \kern -6.5pt
\intop \kern -4pt}\nolimits_{#1}}%
            {\mathop{\vrule width 5pt height
3 pt depth -2.6pt
                    \kern -6pt
\intop \kern -4pt}\nolimits_{#1}}%
            {\mathop{\vrule width 5pt height
3 pt depth -2.6pt
          \kern -6pt \intop \kern -4pt}\nolimits_{#1}}}
          
           \def\bariint_#1{\mathchoice
            {\mathop{\vrule width 15pt
height 3 pt depth -2.5pt
                    \kern -15.8pt
\intop \kern -8pt\intop \kern -4pt}\nolimits_{#1}}%
            {\mathop{\vrule width 9pt height
3 pt depth -2.6pt
                    \kern -10.5pt
\intop \kern -8pt\intop \kern -4pt}\nolimits_{#1}}%
            {\mathop{\vrule width 9pt height
3 pt depth -2.6pt
                    \kern -10pt
\intop \kern -8pt\intop \kern -4pt}\nolimits_{#1}}%
            {\mathop{\vrule width 9pt height
3 pt depth -2.6pt
          \kern -8pt \intop \kern -10pt\intop \kern -4pt}
      \nolimits_{  #1}}}

\def\barintlim_#1{\mathchoice
            {\mathop{\vrule width 6pt
height 3 pt depth -2.5pt
                    \kern -8.8pt
\intop \kern -4pt}\limits_{#1}}%
            {\mathop{\vrule width 5pt height
3 pt depth -2.6pt
                    \kern -6.5pt
\intop \kern -4pt}\limits_{#1}}%
            {\mathop{\vrule width 5pt height
3 pt depth -2.6pt
                    \kern -6pt
\intop \kern -4pt}\limits_{#1}}%
            {\mathop{\vrule width 5pt height
3 pt depth -2.6pt
          \kern -6pt \intop \kern -4pt}\limits_{#1}}}
          
           \def\bariintlim_#1{\mathchoice
            {\mathop{\vrule width 15pt
height 3 pt depth -2.5pt
                    \kern -15.8pt
\intop \kern -8pt\intop \kern -4pt}\limits_{#1}}%
            {\mathop{\vrule width 9pt height
3 pt depth -2.6pt
                    \kern -10.5pt
\intop \kern -8pt\intop \kern -4pt}\limits_{#1}}%
            {\mathop{\vrule width 9pt height
3 pt depth -2.6pt
                    \kern -10pt
\intop \kern -8pt\intop \kern -4pt}\limits_{#1}}%
            {\mathop{\vrule width 9pt height
3 pt depth -2.6pt
          \kern -8pt \intop \kern -10pt\intop \kern -4pt}
      \limits_{  #1}}}
          
\renewcommand{\iint}{\int \kern -8pt\int}       

\numberwithin{equation}{section}
\makeatletter
\@namedef{subjclassname@2020}{\textup{2020} Mathematics Subject Classification}
\makeatother

\title[Complex Gaussian Correlation]{On plurisubharmonic Gaussian correlation of Barthe and Cordero-Erausquin}

\author{Yi C. Huang} 
\address{School of Mathematical Sciences, Nanjing Normal University, Nanjing 210023, People's Republic of China}
\email{Yi.Huang.Analysis@gmail.com}
\urladdr{https://orcid.org/0000-0002-1297-7674}

\date{\today} 

\keywords{Correlation inequality, Gaussian measure, plurisubharmonic functions.}
\subjclass[2020]{Primary 60E15; Secondary 32U05.}  
\thanks{Research of the author is supported by the National NSF grant of China (no. 11801274).
The author would like to thank Sijie Luo (CSU) for helpful communications on concentration inequalities and Fuping Shi (NJNU) for detailed presentations in a weekly Groupe de Travail en Analyse.}

\begin{document}

\begin{abstract}
We simplify the proof of a Gaussian correlation inequality for plurisubharmonic functions found by Barthe and Cordero-Erausquin.
The new observation is a second-order integration-by-parts formula in the complex Gaussian setting.
%
%
%
\end{abstract}

\maketitle


\section{Introduction}

A twice continuously differentiable function $f:\bC^n\rightarrow \bR$ is plurisubharmonic (psh for short) if for all $w\in \bC^n$, the $z$-pointwise non-negativity holds
$$\sum_{1\leq j,k\leq n}\pd^2_{z_j\overline{z_k}}f(z)w_j\overline{w_k}\geq0,$$
where
$$\pd_{z_j}=\frac12(\pd_{x_j}-i\pd_{y_j})\quad\text{and}\quad\pd_{\overline{z_j}}=\frac12(\pd_{x_j}+i\pd_{y_j}).$$ 
Here, $z=x+iy$ with $x,y\in\bRn$.
A function $f$ defined on $\bC^n$ is circular-symmetric if
$$f(e^{i\theta}w)=f(w),\quad \forall\, \theta\in[-\pi,\pi],\quad \,\forall w\in\bC^n.$$
Consider the standard complex Gaussian measure $d\gamma$ on $\bC^n$:
$$d\gamma(w)=\pi^{-n}e^{-|w|^2}d\ell(w),$$
where $d\ell$ denotes the Lebesgue measure on $\bC^n\simeq\bR^{2n}$.

The following remarkable result was recently established in \cite{BarCorEra22}.

\begin{thm}[Barthe and Cordero-Erausquin] \label{thm:BCE}
Let $f, g: \bC^n\rightarrow [-\infty,\infty)$ be two psh functions (with controlled growth at infinity).
If $f$ is circular-symmetric, then
\begin{equation} \label{e:cgc}
\int fgd\gamma\geq \int fd\gamma\int gd\gamma.
\end{equation}
\end{thm}

We refer to Barthe \cite{Bar17} for a nice survey about Gaussian correlation inequality.

\section{Integration by parts}

The proof of \eqref{e:cgc} in \cite{BarCorEra22} is similar to the real case for convex functions
(see Y. Hu \cite{Hu97} and Schmuckenschl\"ager \cite{Sch00}) for the Ornstein-Uhlenbeck operator
$${\bL} f(w):=\frac14\Delta f(w)-\frac12\langle w,\nabla\rangle_{\bC^n} f(w).$$
Using that $\langle w,\nabla\rangle_{\bC^n}=(x+iy)\cdot(\nabla_x-i\nabla_y)$, we can expand $\bL$ as
$${\bL}=\sum_{j=1}^n\frac14\left(\pd^2_{x_jx_j}+\pd^2_{y_jy_j}\right)-\frac12\left(x_j\pd_{x_j}+y_j\pd_{y_j}\right).$$
Since the plurisubharmonicity involves $\pd\bar\pd$-operators, it is then natural to work with
$$L=\sum_{j=1}^n\left(\pd^2_{z_j\overline{z_j}}-\overline{z_j}\pd_{\overline{z_j}}\right)=\sum_{j=1}^ne^{|z|^2}\pd_{{z_j}}\left(e^{-|z|^2}\pd_{\overline{z_j}}\right).$$
For convenience, we introduce the first order vectorial operators
$$\pd_{\overline{z}}=(\pd_{\overline{z_1}},\pd_{\overline{z_2}},\cdots,\pd_{\overline{z_n}})\quad\text{and}\quad \pd_{{z}}=(\pd_{{z_1}},\pd_{{z_2}},\cdots,\pd_{{z_n}}).$$
We have the following integration by parts formulae.

\begin{lem} \label{l:2ibp}
For regular enough functions $f, g:\bC^n\rightarrow \bC$,
$$\int (L f){g}d\gamma=-\int\pd_{\overline{z}}f\cdot{\pd_z g}d\gamma=\int f(\overline{L} g)d\gamma,$$
where $\overline{L}=\sum_{j=1}^n(\pd^2_{z_j\overline{z_j}}-{z_j}\pd_{{z_j}})$.

Moreover, if $f$ is circular-symmetric,
$$\int (L^2 f){g}d\gamma=\sum_{1\leq j,k\leq n}\int\pd^2_{z_j\overline{z_k}}f\cdot{\pd^2_{z_k\overline{z_j}} g}d\gamma.$$
\end{lem}

The first part is standard. For the second part we note the commutativity 
\begin{equation} \label{e:comm}
\pd_{z_j}L=L\pd_{z_j},\quad\forall\,1\leq j\leq n,
\end{equation}
and the following two operator relations
$$L={\bL} +\frac{i}{2}\sum_{j=1}^n\left(y_j\pd_{x_j}-x_j\pd_{y_j}\right)$$
and
$$\overline{L}={\bL} -\frac{i}{2}\sum_{j=1}^n\left(y_j\pd_{x_j}-x_j\pd_{y_j}\right).$$
Thus, $f$ is circular-symmetric implies 
$$Lf=\overline{L}f=\bL f\quad\text{and}\quad L^2f={\bL}^2f=\overline{L}Lf.$$ 
The second part follows by using \eqref{e:comm} and the two identities in the first part.

\begin{rem}
Note that the integration by part formula also implies
$$\int (-L f)\overline{f}d\gamma=\int|\pd_{\overline{z}}f|^2d\gamma\geq0.$$
\end{rem}

\section{A shorter proof of Theorem \ref{thm:BCE}}

Denote by $P_t^{\bL}$ (resp., $P_t^{L}$) the semigroup on $L^2(d\gamma)$ (resp., complex $L^2(d\gamma)$). 
We follow the strategy of Barthe and Cordero-Erausquin \cite{BarCorEra22}. 
For psh functions $f$ and $g$, with $f$ also being circular-symmetric, introduce the $\bR$-valued function
$$\alpha(t):=\int (P_t^{{\bL} }f)gd\gamma,\quad t\geq0.$$
Since $f$ is circular-symmetric, we have
$$\alpha''(t)=\int ({\bL}^2P_t^{{\bL} }f)gd\gamma=\int (L^2P_t^{L}f)gd\gamma.$$ 
Since $P_t^{L}f=P_t^{\bL}f$ is also plurisubharmonic and circular-symmetric, by Lemma \ref{l:2ibp}
$$\begin{aligned}
\alpha''(t)&=\sum_{1\leq j,k\leq n}\int\pd^2_{z_j\overline{z_k}}P_t^{L}f\cdot{\pd^2_{z_k\overline{z_j}} g}d\gamma\\
&=\int\text{Tr}\left(D^2_{\bC}(P_t^{L}f)D^2_{\bC}g\right)d\gamma\geq0.
\end{aligned}$$
Here $D^2_{\bC}h$ denotes the Hermitian matrix $(\pd^2_{z_j\overline{z_k}}h)_{1\leq j,k\leq n}$.
Note that $\alpha$ has a bounded limit $\alpha(\infty)$.
The convexity of $\alpha$ then implies that it is decreasing on $[0,\infty)$ 
and the correlation inequality \eqref{e:cgc} corresponds to $\alpha(0)\geq\alpha(\infty)$.
The proof is complete.

\begin{rem}
In conclusion, our proof requires lighter knowledge about $P_t^{L}$.
\end{rem}

\bibliographystyle{alpha}
 
\bibliography{Hua-ComplexGaussCorr}

\end{document}